# COMPLETE ENUMERATION OF TWO-LEVEL ORTHOGONAL ARRAYS OF STRENGTH $D$ WITH $D + 2$ CONSTRAINTS

By John Stufken and Boxin Tang[1]

*University of Georgia and Simon Fraser University*

Enumerating nonisomorphic orthogonal arrays is an important, yet very difficult, problem. Although orthogonal arrays with a specified set of parameters have been enumerated in a number of cases, general results are extremely rare. In this paper, we provide a complete solution to enumerating nonisomorphic two-level orthogonal arrays of strength $d$ with $d + 2$ constraints for any $d$ and any run size $n = \lambda 2^d$. Our results not only give the number of nonisomorphic orthogonal arrays for given $d$ and $n$, but also provide a systematic way of explicitly constructing these arrays. Our approach to the problem is to make use of the recently developed theory of $J$-characteristics for fractional factorial designs. Besides the general theoretical results, the paper presents some results from applications of the theory to orthogonal arrays of strength two, three and four.

**1. Introduction.** An orthogonal array of size $n$ with $m$ constraints, $s$ levels and strength $d \geq 2$ is an $n \times m$ matrix with entries from a set of $s$ levels, usually taken as $0, 1, \ldots, s - 1$, such that for every $n \times d$ submatrix, each of the $s^d$ level combinations occurs the same number $\lambda$ of times. Such an array is denoted by $OA(n, m, s, d)$. The five parameters $n, m, s, d, \lambda$ cannot vary independently, as the relation $n = \lambda s^d$ must hold. This is why $\lambda$, called the *index* of the array, is suppressed in the notation $OA(n, m, s, d)$.

The importance of orthogonal arrays cannot be overstated. They are intimately related to many other combinatorial objects such as Hadamard matrices, orthogonal Latin squares and error-correcting codes. They are directly useful as fractional factorial designs in factorial experiments arising

Received May 2005; revised February 2006.

[1]Supported by the Natural Sciences and Engineering Research Council of Canada and the NSF.

from scientific and technological investigations. The most fundamental questions concerning orthogonal arrays are when they exist and how to construct them if they do exist. Indeed, these are the two main themes of [13]. Construction and optimality of orthogonal arrays as fractional factorial designs are discussed in [9]. A comprehensive discussion on the use of orthogonal arrays as fractional factorial designs is given by Wu and Hamada [23]. When more than one orthogonal array is available for a given factorial experiment, the experimenter can use a suitable criterion to choose among the candidate arrays. The most commonly used criteria are those of resolution and aberration, along with their generalizations [4, 10, 16, 20, 26] and [24].

This paper considers the enumeration of orthogonal arrays, an important, yet very difficult, problem. For given parameters, finding the complete set of orthogonal arrays is of theoretical importance in its own right. Knowing the complete set of orthogonal arrays allows us to randomly select one to use from all arrays, an important consideration in designing an experiment. In addition, such knowledge of all orthogonal arrays often provides useful experimental data for theorists. Furthermore, having the complete set of arrays available may lead to an unexpected discovery of orthogonal arrays with other useful statistical properties besides those directly offered by the strength of these arrays. Enumerating orthogonal arrays can be greatly simplified by considering only nonisomorphic orthogonal arrays, as all arrays can easily be generated once all nonisomorphic arrays are available. Two orthogonal arrays are said to be *isomorphic* if one can be obtained from the other by a sequence of operations involving permuting rows, columns and levels, and are said to be *nonisomorphic* otherwise.

Despite the importance of the problem, literature on the subject is considerably scarce as compared to that on the existence and construction problems, due to the obvious reason that the enumeration problem is extremely difficult. Here, we survey some main findings. Lam and Tonchev [14] enumerated all $OA(27, 13, 3, 2)$'s. Hedayat, Seiden and Stufken [12] completely classified all $OA(54, 5, 3, 3)$'s. Saturated two-level orthogonal arrays $OA(4\lambda, 4\lambda - 1, 2, 2)$ for $\lambda = 1, 2, \ldots, 6$ are enumerated in [27, 28] and [13], Theorem 7.37. Sun, Li and Ye [21] provided a complete enumeration of $OA(16, m, 2, 2)$'s for all $3 \leq m \leq 15$ and $OA(20, m, 2, 2)$'s for all $3 \leq m \leq 19$. Recently, Bulutoglu and Margot [2] completely enumerated $OA(80, 6, 2, 4)$'s, $OA(96, 7, 2, 4)$'s, $OA(112, 6, 2, 4)$'s and $OA(144, 8, 2, 4)$'s. Chen, Sun and Wu [5] and Xu [25] considered enumeration of linear orthogonal arrays (regular fractional factorial designs) for small run size $n$. Yumiba, Hyodo and Yamamoto [30] classified all $OA(24, 6, 2, 2)$'s that can be obtained from saturated $OA(24, 23, 2, 2)$'s. The above enumeration results are all for specific sets of parameters. The most general result we know of is that of Seiden and Zemach [17], which gives a complete enumeration of two-level orthogonal arrays of strength $d$ with $d+1$ constraints. Two other results of some generality were reported by Fujii, Namikawa and Yamamoto [11], who completely classified $OA(2^{d+1}, d+2, 2, d)$'s and $OA(2^{d+2}, d+3, 2, d)$'s. Orthogonal arrays are special cases of balanced arrays. Srivastava [19] and Shirakura [18] presented the numbers of nonisomorphic balanced arrays of strength $d$ for $m = d+1, d+2$ and $d+3$ factors. However, their definition of isomorphism is different, in that only row and column permutations are used. In fact, it would not make sense to consider symbol permutations in defining isomorphic balanced arrays, as permuting the symbols of a column in a balanced array may result in an unbalanced array.

In this paper, we provide a complete solution to the problem of enumerating two-level orthogonal arrays of strength $d$ with $d+2$ constraints. The results are general in that they are for any strength $d$ and any index $\lambda$. Our approach is to make use of the recently developed theory of $J$-characteristics for fractional factorial designs [22]. The remainder of the article is structured as follows. Section 2 reviews $J$-characteristics, discusses their application to orthogonal arrays and characterizes nonisomorphic arrays using $J$-characteristics. Section 3 presents our main results through four theorems. Applications to orthogonal arrays of strength two, three and four are also discussed here. Section 4 concludes the paper with a summary of our findings and a discussion of possible future work.

**2. Enumeration method.** Our method of enumerating two-level orthogonal arrays is based on the recently developed theory of $J$-characteristics for fractional factorial designs. Section 2.1 briefly reviews fractional factorial designs and their $J$-characteristics, and presents a key identity linking a fractional factorial with its $J$-characteristics. Section 2.2 specializes $J$-characteristics for orthogonal arrays of strength $d$ with $d+2$ constraints. Section 2.3 considers how to characterize nonisomorphic orthogonal arrays of strength $d$ with $d+2$ constraints using $J$-characteristics, paving the way for the main results to come in Section 3.

2.1. *Fractional factorial designs and their J-characteristics.* For any $s \subseteq Z_m = \{1, \ldots, m\}$, define a row vector

(1) $$\mathbf{r}_s = (r_{s1}, \ldots, r_{sm}),$$

where $r_{sj} = -1$ if $j \in s$ and $r_{sj} = +1$ otherwise. A complete $2^m$ factorial can then be represented by the $2^m \times m$ matrix

(2) $$\mathbf{C} = (\mathbf{r}_\phi^T, \mathbf{r}_1^T, \mathbf{r}_2^T, \mathbf{r}_{12}^T, \mathbf{r}_3^T, \mathbf{r}_{13}^T, \mathbf{r}_{23}^T, \mathbf{r}_{123}^T, \mathbf{r}_4^T, \mathbf{r}_{14}^T, \ldots)^T,$$

where, for simplicity, we use $\mathbf{r}_1$ for $\mathbf{r}_{\{1\}}$, $\mathbf{r}_{12}$ for $\mathbf{r}_{\{1,2\}}$ and so on. Let $\mathbf{h}_j$ denote the $j$th column of $\mathbf{C}$ in (2), $j = 1, \ldots, m$, that is,

$$\mathbf{C} = (\mathbf{h}_1, \ldots, \mathbf{h}_m).$$

The common component of $\mathbf{h}_j$ and $\mathbf{r}_s$ is $r_{sj}$, as in (1). By the $s$th row of $\mathbf{C}$, we mean $\mathbf{r}_s$ and, similarly, by the $s$th component of $\mathbf{h}_j$, we mean $r_{sj}$. For any $t \subseteq Z_m$, the Hadamard product $\mathbf{h}_t$ of the columns $\mathbf{h}_j$ with $j \in t$ is a column vector of length $2^m$, the $s$th component of which is $\prod_{j \in t} r_{sj}$. Let $\mathbf{H}$ be the Hadamard matrix of order $2^m$ formed by collecting in Yates order all of the Hadamard products of the columns of $\mathbf{C}$, that is,

$$\mathbf{H} = (\mathbf{h}_\phi, \mathbf{h}_1, \mathbf{h}_2, \mathbf{h}_{12}, \mathbf{h}_3, \mathbf{h}_{13}, \mathbf{h}_{23}, \mathbf{h}_{123}, \mathbf{h}_4, \mathbf{h}_{14}, \ldots). \tag{3}$$

Deleting $\mathbf{h}_\phi$ from $\mathbf{H}$ gives a saturated factorial in $m$ independent columns.

A two-level fractional factorial design of $n$ runs for $m$ factors is an $n \times m$ matrix $\mathbf{D} = (d_{ij})$ with $d_{ij} = \pm 1$, where each row of $\mathbf{D}$ corresponds to a run and each column to a factor. A row permutation of $\mathbf{D}$ gives the same design. Let $N_s$ denote the number of times that a run $\mathbf{r}_s$ in (1) occurs in $\mathbf{D}$. This very $N_s$ is called the *indicator function* by some authors; the reader is referred to [29] for a more extensive discussion. Design $\mathbf{D}$ can then be equivalently described by the vector of length $2^m$,

$$\mathbf{N} = (N_\phi, N_1, N_2, N_{12}, N_3, N_{13}, N_{23}, N_{123}, N_4, N_{14}, \ldots)^{\mathrm{T}}, \tag{4}$$

where, for example, $N_{12}$ is shorthand notation for $N_{\{1,2\}}$. As a matrix, $\mathbf{D}$ is determined by $\mathbf{N}$ up to a row permutation and, as a design, $\mathbf{D}$ is completely determined by $\mathbf{N}$. Conversely, any vector $\mathbf{N}$ of length $2^m$ with nonnegative integers as its components defines a design with $n = \sum_{s \subseteq Z_m} N_s$ runs.

For any $t \subseteq Z_m$, let

$$J_t = \sum_{i=1}^n \prod_{j \in t} d_{ij} = \sum_{s \subseteq Z_m} h_{st} N_s, \tag{5}$$

where $\mathbf{H} = (h_{st})$ is the Hadamard matrix in (3) and $h_{st}$ denotes the $s$th entry of $\mathbf{h}_t$, the $t$th column of $\mathbf{H}$. The $J_t$ values defined in (5) for all subsets $t$ of $Z_m$ are called the *J-characteristics* of design $\mathbf{D}$. Let

$$\mathbf{J} = (J_\phi, J_1, J_2, J_{12}, J_3, J_{13}, J_{23}, J_{123}, J_4, J_{14}, \ldots)^{\mathrm{T}}. \tag{6}$$

The set of $J$-characteristics as in (6) completely determines a general factorial, just as the defining relation does for a regular factorial. This is given in the following lemma.

LEMMA 1. $\mathbf{N} = 2^{-m} \mathbf{H} \mathbf{J}$, *that is,*

$$N_s = 2^{-m} \sum_{t \subseteq Z_m} h_{st} J_t, \tag{7}$$

*where* $\mathbf{H} = (h_{st})$ *is as given in (3).*

This result is given in [22] and plays a key role in this paper. As discussed in [22], $J$-characteristics directly capture the projection properties of a factorial design and, because of this, they lead to successful generalizations of minimum aberration from regular factorials to general factorials [7]. Further work on generalized minimum aberration can be found in [26] and [3]. Finally, we note that, mathematically, **J** is simply the Hadamard transform of **N**.

2.2. *J-characteristics for orthogonal arrays.* For two-level orthogonal arrays, it is convenient to use $\pm 1$ to denote the two levels. A general two-level fractional factorial design $\mathbf{D} = (d_{ij})_{n \times m}$, where $d_{ij} = \pm 1$ as introduced in Section 2.1, may or may not be an orthogonal array of strength $d \geq 2$, but, if it is, then its $J$-characteristics have some special properties. The next three lemmas are devoted to these properties of $J$-characteristics of an orthogonal array.

LEMMA 2. *Design* **D** *is an orthogonal array of strength* $d$ *if and only if* $J_t = 0$ *for all* $t \subseteq Z_m$ *with* $|t| \leq d$, *where* $t \neq \phi$.

Lemma 2 can easily be verified using (5) and (7) and was given earlier in [7].

LEMMA 3. *Suppose that* **D** *is an* $OA(n, m, 2, d)$, *where* $n = \lambda 2^d$ *and* $m \geq d + 2$. *Then:*

(i) *for any* $t \subseteq Z_m$, *it holds that* $J_t = \mu_t 2^d$ *for some integer* $\mu_t$;
(ii) *if* $\lambda$ *is even, then* $\mu_t$ *is even;*
(iii) *if* $\lambda$ *is odd and* $d$ *is even, then* $\mu_t$ *is odd for* $|t| = d+1, d+2$;
(iv) *if* $\lambda$ *is odd and* $d$ *is odd, then* $\mu_t$ *is odd for* $|t| = d+1$ *and even for* $|t| = d+2$.

PROOF. From (1), we see that $\mathbf{r}_\phi = (1, \ldots, 1)$, implying that $h_{\phi t} = 1$ for any $t$. Taking $s = \phi$ in (7) gives

$$N_\phi = 2^{-m} \sum_{t \subseteq Z_m} J_t. \tag{8}$$

For given $w \subseteq Z_m$, consider the projection design $\mathbf{D}_w$ consisting of the columns of **D** in $w$. Now applying (8) to this projection design $\mathbf{D}_w$, we obtain $N_\phi = 2^{-|w|} \sum_{t \subseteq w} J_t$. From Lemma 2, we have $J_t = 0$ for all $t$ with $1 \leq |t| \leq d$. Further, note that $J_\phi = n = \lambda 2^d$. Therefore,

$$J_w = 2^{|w|} N_\phi - \lambda 2^d - \sum_{t \subset w, d+1 \leq |t| \leq |w|-1} J_t. \tag{9}$$

Based on (9), Lemma 3 can easily be established. □

For the special case $d = 2$, results similar to those in Lemma 3 were provided in [8]. Lemmas 2 and 3 discuss the possible values of $J_t$ for a given $t$. The following result states that the values of $J_t$ for two different subsets $t_1$ and $t_2$ are, in fact, related, provided that the two subsets are close to each other.

LEMMA 4. *Let $\mathbf{D}$ be an $OA(n, m, 2, d)$. We must have*

$$|J_{t_1}| + |J_{t_2}| \leq n,$$

*provided that $1 \leq |t_1 \triangle t_2| \leq d$, where $t_1 \triangle t_2 = (t_1 \setminus t_2) \cup (t_2 \setminus t_1)$.*

PROOF. Let $\mathbf{d}_j$ be the $j$th column of $\mathbf{D}$. Let $\mathbf{d}_{t_1}$ be the Hadamard product of the columns $\mathbf{d}_j$ for $j \in t_1$. Define $\mathbf{d}_{t_2}$ similarly. Consider the design $\mathbf{D}^*$ given by the two column vectors $\mathbf{d}_{t_1}$ and $\mathbf{d}_{t_2}$, that is, $\mathbf{D}^* = (\mathbf{d}_{t_1}, \mathbf{d}_{t_2})$. Obviously, the $J$-characteristics for column vectors $\mathbf{d}_{t_1}$ and $\mathbf{d}_{t_2}$ are $J_{t_1}$ and $J_{t_2}$, respectively, and the $J$-characteristic for the two combined is $J_{t_1 \triangle t_2}$. Because $1 \leq |t_1 \triangle t_2| \leq d$ and $\mathbf{D}$ is an orthogonal array of strength $d$, we must have $J_{t_1 \triangle t_2} = 0$. We apply Lemma 1 to design $\mathbf{D}^*$. Note that the corresponding Hadamard matrix $\mathbf{H}$ in (3) is now given by

$$\mathbf{H} = [(1,1,1,1)^{\mathrm{T}}, (1,-1,1,-1)^{\mathrm{T}}, (1,1,-1,-1)^{\mathrm{T}}, (1,-1,-1,1)^{\mathrm{T}}]^{\mathrm{T}}.$$

Noting that $N_s$ in Lemma 1 is nonnegative, we obtain

$$n + J_{t_1} + J_{t_2} \geq 0, \qquad n - J_{t_1} + J_{t_2} \geq 0,$$
$$n + J_{t_1} - J_{t_2} \geq 0, \qquad n - J_{t_1} - J_{t_2} \geq 0.$$

This shows that $|J_{t_1}| + |J_{t_2}| \leq n$. □

The seemingly unremarkable result in Lemma 4 is actually very powerful. For example, a key result in [6], Lemma 2.2, is a direct consequence of our Lemma 4. In what follows, we document an interesting application of Lemma 4 to orthogonal arrays of strength 2, the most useful situation in the statistical design of experiments.

COROLLARY 1. *Let $\mathbf{D}$ be an $OA(4\lambda, m, 2, 2)$ with odd $\lambda$. Then we must have $|J_t| \leq 4\lambda - 4$ for all $t$ with $1 \leq |t| \leq m - 1$. For $t = Z_m$, the possibility that $|J_t| = 4\lambda$ can occur only when the two conditions simultaneously hold that (i) $m \equiv 3 \pmod{4}$ and (ii) $\mathbf{D}$ cannot be embedded into an $OA(4\lambda, m + 1, 2, 2)$.*

We omit the proof, which involves using a result in [8], Proposition 1, and is otherwise straightforward. If $|J_t| = 4\lambda$, then the set of columns $\mathbf{d}_j$ with $j \in t$ is said to form a *defining word*, which results in full effect aliasing. Corollary 1 effectively states that when the run size $n = 4\lambda$ is not a multiple of 8, then, except in rare cases, an orthogonal array of strength 2 does not have a defining word.

2.3. *Characterizing nonisomorphic $OA(n, d+2, 2, d)$'s using $J$-characteristics.* As discussed in Section 2.1, a fractional factorial design $\mathbf{D}$ is completely defined, up to a row permutation, by the vector $\mathbf{N}$ in (4). Lemma 1 asserts that $\mathbf{N}$ is uniquely determined by $\mathbf{J}$ in (6), and vice versa. Therefore, up to a row permutation, a factorial design $\mathbf{D}$ is completely characterized by its $J$-characteristics. When design $\mathbf{D}$ is an $OA(n, d+2, 2, d)$, Lemma 2 asserts that $J_\phi = n$, and $J_t = 0$ for all nonempty $t$ with $|t| \leq d$. So, the vector $\mathbf{J}$ in (6) can be equivalently described by

$$\mathbf{J} = (J_{t_1}, \ldots, J_{t_{m+1}}), \tag{10}$$

where $t_j = Z_m \setminus \{m+1-j\}$ for $j = 1, \ldots, m$ and $t_{m+1} = Z_m$. For convenience, we use the same notation $\mathbf{J}$ for this shortened version of the original vector $\mathbf{J}$. No confusion will arise, as the new vector $\mathbf{J}$ in (10) is always used in the context of $OA(n, d+2, 2, d)$'s. We call the vector $\mathbf{J}$ in (10) the $\mathbf{J}$-*vector* of the $OA(n, d+2, 2, d)$ under consideration. Two orthogonal arrays are said to be *statistically equivalent* if one can be obtained from the other by a row permutation [13]. Thus, two statistically equivalent $OA(n, d+2, 2, d)$'s have the same $\mathbf{J}$-vector and two $OA(n, d+2, 2, d)$'s that are not statistically equivalent have different $\mathbf{J}$-vectors. The preceding discussion leads to the following conclusion.

LEMMA 5. *There is a one-to-one correspondence between the classes of statistically equivalent $OA(n, d+2, 2, d)$'s and their $\mathbf{J}$-vectors.*

It is tempting to use the $\mathbf{J}$-vector to identify nonisomorphic arrays, but in its present form, the $\mathbf{J}$-vector does not serve this purpose. This is because two isomorphic arrays will have two different $\mathbf{J}$-vectors if they are not statistically equivalent; for a given array, both permuting its columns and switching the signs of its columns affect its $\mathbf{J}$-vector. Lemma 5 has taken care of isomorphism due to row-permuting. We now consider isomorphism due to column-permuting and sign-switching. For a given array, there are $m!$ ways of permuting its columns. Let $\delta_j$ be the indicator of whether the $j$th column is sign-switched, with $\delta_j = -1$ denoting that there is sign-switching and $\delta_j = +1$ denoting that there is no sign-switching. A vector $\boldsymbol{\delta} = (\delta_1, \ldots, \delta_m)$ then gives one way of sign-switching the columns of the array. Thus, the total number of ways of sign-switching the columns is $2^m$. Starting from a given array, column-permuting and sign-switching together generate $m!2^m$ arrays in all, some of which may be statistically equivalent. From the **J**-vectors of these $m!2^m$ arrays, we will choose a unique **J**-vector, and therefore a unique array, to represent the whole class of isomorphic arrays.

We first present a lemma which will be needed in what follows and also in Section 3. Consider the fractional factorial design **B** obtained by collecting columns $\mathbf{b}_1, \ldots, \mathbf{b}_{m+1}$ from the Hadamard matrix **H** in (3), that is,

$$\mathbf{B} = (\mathbf{b}_1, \ldots, \mathbf{b}_{m+1}), \tag{11}$$

where $\mathbf{b}_j = \mathbf{h}_{t_j}$ with $t_j = Z_m \setminus \{m+1-j\}$ for $j = 1, \ldots, m$ and $t_{m+1} = Z_m$. Recall that a defining word in a fractional factorial design is a set of columns whose Hadamard product is **I**, a column of all $+1$'s.

LEMMA 6. *Design **B** has exactly one defining word, given by $\mathbf{b}_1 \cdots \mathbf{b}_{m+1} = \mathbf{I}$ when $m$ is even and $\mathbf{b}_1 \cdots \mathbf{b}_m = \mathbf{I}$ when $m$ is odd. Therefore, for even $m$, deleting any column from **B** gives a full factorial design and for odd $m$, deleting any column except for $\mathbf{b}_{m+1}$ also gives a full factorial.*

The proof is simple and is thus omitted.

We now return to orthogonal arrays. Let $\mathbf{D} = (\mathbf{d}_1, \ldots, \mathbf{d}_m)$ be an $OA(n, m, 2, d)$ where $m = d + 2$ and $\mathbf{d}_j$ denotes the $j$th column of **D**. As discussed earlier, a vector $\boldsymbol{\delta} = (\delta_1, \ldots, \delta_m)$ with $\delta_j = \pm 1$ gives one way of switching the signs of the columns. Sign-switching **D** by $\boldsymbol{\delta} = (\delta_1, \ldots, \delta_m)$ results in a new array $\mathbf{D}_{\boldsymbol{\delta}} = (\delta_1 \mathbf{d}_1, \ldots, \delta_m \mathbf{d}_m)$. Obviously, we have $J_{t_j}(\mathbf{D}_{\boldsymbol{\delta}}) = \delta_{t_j} J_{t_j}(\mathbf{D})$ for $j = 1, \ldots, m+1$, where $\delta_{t_j} = \prod_{i \in t_j} \delta_i$. Since the $2^m$ vectors $\boldsymbol{\delta} = (\delta_1, \ldots, \delta_m)$ collectively form a full factorial, the resulting $2^m$ vectors $(\delta_{t_1}, \ldots, \delta_{t_{m+1}})$ are exactly the $2^m$ runs of design **B** given in (11).

Now, suppose $m$ is even. The situation involving odd $m$ will be dealt with later.

Case (i): If $|J_{t_{m+1}}| \leq \min_{1 \leq j \leq m} |J_{t_j}|$, select a $\boldsymbol{\delta} = (\delta_1, \ldots, \delta_m)$ such that $\delta_{t_j} = -1$ if $J_{t_j} > 0$ and $\delta_{t_j} = +1$ if $J_{t_j} \leq 0$, for $j = 1, \ldots, m$. This is possible because for even $m$, the $2^m$ vectors $(\delta_{t_1}, \ldots, \delta_{t_m})$ form a full factorial according to Lemma 6. Such a choice of $\boldsymbol{\delta} = (\delta_1, \ldots, \delta_m)$ leads to an array, the **J**-vector of which satisfies

$$J_{t_j} \leq 0 \qquad \text{for } j = 1, \ldots, m. \tag{12}$$

Clearly, permuting the columns of an array results in a permutation of $J_{t_1}, \ldots, J_{t_m}$. It is also easy to see that any permutation of $J_{t_1}, \ldots, J_{t_m}$ corresponds to a permutation of the columns. For the array satisfying (12), we can permute its columns to obtain an array whose **J**-vector satisfies

$$J_{t_1} \leq \cdots \leq J_{t_m} \leq 0. \tag{13}$$

The **J**-vector resulting from the above procedure is unique, as is its corresponding array. We note that although there may be more than one $\boldsymbol{\delta}$-vector or more than one permutation of columns that can lead to a **J**-vector satisfying (13), this **J**-vector is always unique. For example, if $J_{t_1} = 0$, we can also use a $\boldsymbol{\delta}$-vector having $\delta_{t_1} = -1$ instead of $\delta_{t_1} = +1$, which will not have any effect on the resulting **J**-vector since both $J_{t_1}$ and $J_{t_{m+1}}$ are still zero.

Case (ii): If $|J_{t_{m+1}}| > \min_{1 \leq j \leq m} |J_{t_j}| = |J_{t_{j'}}|$, select a $\boldsymbol{\delta} = (\delta_1, \ldots, \delta_m)$ such that $\delta_{t_j} = -1$ if $J_{t_j} > 0$ and $\delta_{t_j} = +1$ if $J_{t_j} \leq 0$ for $j = 1, \ldots, j'-1, j'+1, \ldots, m+1$. This is possible because for even $m$, the $2^m$ vectors $(\delta_{t_1}, \ldots, \delta_{t_{j'-1}}, \delta_{t_{j'+1}}, \ldots, \delta_{t_{m+1}})$ also form a full factorial according to Lemma 6. This choice of $\boldsymbol{\delta}$ gives an array whose **J**-vector satisfies

$$J_{t_{m+1}} < 0 \quad \text{and} \quad J_{t_j} \leq 0 \qquad \text{for } j = 1, \ldots, j'-1, j'+1, \ldots, m.$$

Now, permute $J_{t_1}, \ldots, J_{t_m}$ by first moving $J_{t_{j'}}$ to the end of the sequence and then arranging the others in ascending order. We now obtain an array whose **J**-vector satisfies

(14) $$J_{t_1} \leq \cdots \leq J_{t_{m-1}} \leq 0 \quad \text{and} \quad J_{t_{m+1}} < 0.$$

Once again, the resulting **J**-vector is unique, as is the corresponding array. It is possible that more than one $j'$ satisfies $\min_{1 \leq j \leq m} |J_{t_j}| = |J_{t_{j'}}|$. We can easily check that it has no effect on the final outcome.

The above results can be summarized in the following proposition.

PROPOSITION 1. *When $m$ is even, every class of isomorphic orthogonal arrays of strength $d$ with $m = d+2$ constraints contains a unique array whose **J**-vector satisfies either the condition in (15) or that in (16) as given by*

(15) $$J_{t_1} \leq \cdots \leq J_{t_m} \leq -|J_{t_{m+1}}|,$$

(16) $$J_{t_1} \leq \cdots \leq J_{t_{m-1}} \leq -|J_{t_m}|, \qquad J_{t_{m+1}} < -|J_{t_m}|.$$

Proposition 1 is immediate as (15) is (13) under case (i) and (16) is (14) under case (ii).

The situation for odd $m$ is considerably simpler. Let $j'$ be such that $|J_{t_{j'}}| = \min_{1 \leq j \leq m} |J_{t_j}|$. We select a $\boldsymbol{\delta} = (\delta_1, \ldots, \delta_m)$ such that $\delta_{t_j} = -1$ if $J_{t_j} > 0$ and $\delta_{t_j} = +1$ if $J_{t_j} \leq 0$ for $j = 1, \ldots, j'-1, j'+1, \ldots, m+1$. This is possible because for odd $m$, the $2^m$ vectors $(\delta_{t_1}, \ldots, \delta_{t_{j'-1}}, \delta_{t_{j'+1}}, \ldots, \delta_{t_{m+1}})$ form a full factorial according to Lemma 6. This choice of $\boldsymbol{\delta}$ gives an array whose **J**-vector satisfies

$$J_{t_{m+1}} \leq 0 \quad \text{and} \quad J_{t_j} \leq 0 \qquad \text{for } j = 1, \ldots, j'-1, j'+1, \ldots, m.$$

Now, permute $J_{t_1}, \ldots, J_{t_m}$ by first moving $J_{t_{j'}}$ to the end of the sequence and then arranging the others in ascending order. We now obtain an array, the **J**-vector of which satisfies

$$J_{t_1} \leq \cdots \leq J_{t_{m-1}} \leq 0 \quad \text{and} \quad J_{t_{m+1}} \leq 0.$$

Using the same arguments as before, we conclude that the resulting **J**-vector is unique.

PROPOSITION 2. *When $m$ is odd, every class of isomorphic orthogonal arrays of strength $d$ and with $m = d + 2$ constraints contains a unique array whose **J**-vector satisfies*

$$\text{(17)} \qquad J_{t_1} \leq \cdots \leq J_{t_{m-1}} \leq -|J_{t_m}|, \qquad J_{t_{m+1}} \leq 0.$$

To enumerate nonisomorphic $OA(n, d+2, 2, d)$'s, Propositions 1 and 2 conclude that we only need to consider $OA(n, d+2, 2, d)$'s whose **J**-vectors satisfy (15) or (16) for even $m$ and (17) for odd $m$. For convenience, a **J**-vector is called a **J**$^*$-*vector* if it satisfies (15) or (16) for even $m$ and (17) for odd $m$. Since a **J**$^*$-vector uniquely determines a class of isomorphic $OA(n, d+2, 2, d)$'s, two $OA(n, d+2, 2, d)$'s are nonisomorphic if their **J**$^*$-vectors are different. We therefore establish a one-to-one correspondence between nonisomorphic $OA(n, d+2, 2, d)$'s and their **J**$^*$-vectors. Enumerating orthogonal arrays of strength $d$ with $m = d + 2$ constraints is equivalent to enumerating their **J**$^*$-vectors.

**3. Enumeration results.** Our results on enumerating $OA(\lambda 2^d, d+2, 2, d)$'s are presented in Sections 3.1, 3.2, 3.3 and 3.4, respectively, dealing with the four cases (i) even $d$ and odd $\lambda$, (ii) even $d$ and even $\lambda$, (iii) odd $d$ and odd $\lambda$ and (iv) odd $d$ and even $\lambda$. Since these four sections are structurally similar, we shall only preview Section 3.1 here. The main result is contained in Theorem 1, which establishes a one-to-one correspondence between the set of nonisomorphic $OA(\lambda 2^d, d+2, 2, d)$'s with even $d$ and odd $\lambda$ and the set of solutions to equation (18) with its unknowns satisfying certain conditions. Example 1 and Corollary 2 give a first taste of the power of the theorem. Through Lemmas 7 and 8, we then provide a systematic method for explicitly obtaining all solutions to equation (18). Finally, some results from the application of Theorem 1 and Lemmas 7 and 8 are presented in Tables 1 and 2.

3.1. *Enumerating $OA(\lambda 2^d, d+2, 2, d)$'s for even $d$ and odd $\lambda$.* Recall that $m = d + 2$. Consider the equation

$$\text{(18)} \qquad \lambda + u_1 + \cdots + u_{m+1} = 4k,$$

where $k$ and the $u_j$'s are the unknowns, $k \geq 0$ is an integer and the $u_j$'s are odd integers with $|u_j| \leq \lambda - 2$ for all $j = 1, \ldots, m+1$. Further, the $u_j$'s must satisfy one of the conditions

$$\text{(19)} \qquad u_1 \leq \cdots \leq u_m \leq -|u_{m+1}|,$$

$$\text{(20)} \qquad u_1 \leq \cdots \leq u_{m-1} \leq -|u_m|, \qquad u_{m+1} \leq -|u_m| - 2.$$

It is clear that (19) and (20) are mutually exclusive. Note that both the solutions to (18) under (19) and those to (18) under (20) are considered in the following theorem.

THEOREM 1. *Suppose that $d$ is even and $\lambda$ is odd. We then have that:*

(i) *each solution $(u_1, \ldots, u_{m+1}, k)$ to equation (18) under either (19) or (20) determines an $OA(\lambda 2^d, d+2, 2, d)$ with $\mathbf{J}^*$-vector given by $J_{t_j} = 2^d u_j$ for $j = 1, \ldots, m+1$;*

(ii) *the complete set of nonisomorphic $OA(\lambda 2^d, d+2, 2, d)$'s is given by collecting the arrays obtained in (i) over all the solutions to equation (18).*

PROOF. (i) Let $(u_1, \ldots, u_{m+1}, k)$ be a solution to (18) with all of the accompanying conditions satisfied. Consider

$$(21) \quad \mathbf{J}^* = (J_{t_1}, \ldots, J_{t_{m+1}}), \qquad \text{where } J_{t_j} = 2^d u_j \text{ for } j = 1, \ldots, m+1.$$

We will show that the vector $\mathbf{J}^* = (J_{t_1}, \ldots, J_{t_{m+1}})$ in (21) is indeed the $\mathbf{J}^*$-vector of an $OA(\lambda 2^d, d+2, 2, d)$. Obviously, (19) implies (15) and (20) implies (16). From Lemma 1, it is clear that our goal is achieved if we can show that

$$N_s = 2^{-m}\left(n + \sum_{j=1}^{m+1} h_{st_j} J_{t_j}\right)$$

is a nonnegative integer for all $s \subseteq Z_m$. Noting that $m = d+2$ and $n = \lambda 2^d$, we have

$$N_s = \left(\lambda + \sum_{j=1}^{m+1} h_{st_j} u_j\right)\bigg/4 = \left(\lambda + \sum_{j=1}^{m+1} u_j\right)\bigg/4 + \sum_{j=1}^{m+1} [h_{st_j} - 1]u_j/4.$$

Since $\lambda + \sum_{j=1}^{m+1} u_j = 4k$ and $h_{st_j} = \pm 1$, we obtain

$$N_s = k - \sum_{j: h_{st_j} = -1} u_j/2 = k - \Delta/2,$$

where $\Delta = \sum_{j: h_{st_j} = -1} u_j$. Since $k \geq 0$ is an integer, if we can show that $\Delta$ is a nonpositive even integer, then $N_s$ must be a nonnegative integer. From Lemma 6, we see that the number of $j$'s with $h_{st_j} = -1$ has to be even because $m$ is even and $(h_{st_1}, \ldots, h_{st_{m+1}})$ is a run of design $\mathbf{B}$ in (11). We now obtain that $\Delta$ is even, since all $u_j$'s are odd and so $\Delta = \sum_{j: h_{st_j} = -1} u_j$ is a sum of an even number of odd numbers. For $s = \phi$, we have $N_s = k \geq 0$. For a nonempty $s$, the number of $j$'s with $h_{st_j} = -1$ is a positive even integer, implying that $\sum_{j: h_{st_j} = -1} u_j$ is a sum of at least two $u_j$'s. From (19) or (20), we see that the only $u_j$ that can possibly be positive has the smallest absolute value among all $u_j$'s. [This possibly positive $u_j$ is $u_{m+1}$ if (19) holds and is $u_m$ if (20) holds.] We can now conclude that $\Delta \leq 0$.

(ii) Part (i) of Theorem 1 states that a solution $(u_1, \ldots, u_{m+1}, k)$ to (18) determines an $OA(\lambda 2^d, d+2, 2, d)$ with $\mathbf{J}^*$-vector given by $J_{t_j} = 2^d u_j$. Clearly, two different solutions to (18) give two different $\mathbf{J}^*$-vectors, representing two nonisomorphic $OA(\lambda 2^d, d+2, 2, d)$'s from Proposition 1. Part (ii) of Theorem 1 is established if we can show that every $OA(\lambda 2^d, d+2, 2, d)$ can be obtained from a solution to (18). Now, suppose an $OA(\lambda 2^d, d+2, 2, d)$ has a $\mathbf{J}^*$-vector given by $\mathbf{J}^* = (J_{t_1}, \ldots, J_{t_{m+1}})$. Let $u_j = J_{t_j}/2^d$ for $j = 1, \ldots, m+1$. From Lemma 3(iii), we know that the $u_j$'s are odd integers. Lemma 4 further implies that $|u_j| + |u_{j'}| \leq \lambda$ for $j \neq j'$. Since $u_j$, $u_{j'}$ and $\lambda$ are all odd, we must have $|u_j| \leq \lambda - 2$ for all $j = 1, \ldots, m+1$. Now, consider

$$N_\phi = 2^{-m}\left(n + \sum_{j=1}^{m+1} J_{t_j}\right) = (\lambda + u_1 + \cdots + u_{m+1})/4.$$

We then see that $(u_1, \ldots, u_{m+1}, k)$, where $k = N_\phi$, is a solution to equation (18). The condition in (19) or (20) is met due to Proposition 1. $\square$

For even $d$ and odd $\lambda$, Theorem 1 establishes a one-to-one correspondence between the set of nonisomorphic $OA(\lambda 2^d, d+2, 2, d)$'s and the set of solutions to equation (18) subject to the given conditions. From the above proof, we see that $N_s \geq N_\phi = k \geq 0$, so the $k$ in a solution $(u_1, \ldots, u_{m+1}, k)$ to (18) has a simple interpretation: the corresponding array contains $k$ copies of a full factorial, with $k$ being the largest integer for which such a statement can be made.

EXAMPLE 1. First, consider the case $\lambda = 3$, $d = 2$, $n = 12$ and $m = 4$. Equation (18) has a unique solution under (19), given by $(u_1, u_2, u_3, u_4, u_5, k) = (-1, -1, -1, -1, 1, 0)$, and has no solution under (20), so the $OA(12, 4, 2, 2)$ is unique up to isomorphism. For the case $\lambda = 5$, $d = 2$, $n = 20$ and $m = 4$, equation (18) has two solutions $(-1, -1, -1, -1, -1, 0)$, $(-3, -1, -1, -1, 1, 0)$ under (19) and one solution $(-1, -1, -1, 1, -3, 0)$ under (20). Therefore, there are three nonisomorphic $OA(20, 4, 2, 2)$'s. For $\lambda = 7$, $d = 2$, $n = 28$ and $m = 4$, there are seven nonisomorphic $OA(28, 4, 2, 2)$'s given by the solutions

$$\begin{aligned}(u_1, &u_2, u_3, u_4, u_5, k) \\ &= (-1, -1, -1, -1, 1, 1), (-3, -3, -1, -1, 1, 0), \\ &\quad (-5, -1, -1, -1, 1, 0), (-3, -1, -1, -1, -1, 0), \\ &\quad (-1, -1, -1, -1, -3, 0), (-3, -1, -1, 1, -3, 0), \\ &\quad (-1, -1, -1, 1, -5, 0),\end{aligned}$$

where the first four are solutions under (19) and the remaining three are solutions under (20). Similarly, we have the following conclusion for orthogonal arrays of strength 4. There exist no $OA(48, 6, 2, 4)$, one nonisomorphic $OA(80, 6, 2, 4)$, three nonisomorphic $OA(112, 6, 2, 4)$'s and seven nonisomorphic $OA(144, 6, 2, 4)$'s. In fact, the following general result holds.

COROLLARY 2. *Let $g(\lambda, d)$ be the number of nonisomorphic $OA(\lambda 2^d, d + 2, 2, d)$'s with even $d$ and odd $\lambda$. Then $g(\lambda, d) = 0$ for $\lambda \leq d - 1$ and $g(\lambda, d) = 1, 3, 7$ for $\lambda = d + 1, d + 3, d + 5$, respectively.*

For odd $d$ and odd $\lambda$, a conclusion similar to Corollary 2 also holds; see Corollary 4 in Section 3.3. Nonexistence of an $OA(\lambda 2^d, d + 2, 2, d)$ for $\lambda \leq d - 1$ has previously been established by Blum, Schatz and Seiden [1]; see also Theorem 2.29 in [13]. Corollary 2 can easily be verified directly using Theorem 1 and can also be obtained from the following Lemmas 7 and 8.

The complete set of solutions to (18) under either (19) or (20) can, in fact, be obtained explicitly in a systematic fashion. For ease of presentation, we first introduce some notation. Let $Z[a, b]$ denote the set of integers $x$ such that $a \leq x \leq b$ and $O[a, b]$ denote the set of odd integers $x$ such that $a \leq x \leq b$. Let $S_1$ be the set of solutions to (18) under (19) and $S_2$ be the set of solutions to (18) under (20). We then have the following results.

LEMMA 7. *For even $d$ and odd $\lambda$, if $\lambda \leq d - 1$, then equation (18) has no solution under (19). If $\lambda \geq d + 1$, then equation (18) has at least one solution under (19) and the complete set $S_1$ of solutions is given by*

$$k \in Z[0, (\lambda - d - 1)/4],$$
$$u_{m+1} \in O[-(\lambda - 4k)/(m + 1), (\lambda - 4k)/(m - 1)],$$
$$u_m \in O[-(\lambda - 4k + u_{m+1})/m, -|u_{m+1}|],$$
$$u_j \in O[-(\lambda - 4k + u_{j+1} + \cdots + u_{m+1})/j, u_{j+1}]$$
$$\text{for } j = m - 1, m - 2, \ldots, 2,$$
$$u_1 = -(\lambda - 4k + u_2 + \cdots + u_{m+1}).$$

LEMMA 8. *For even $d$ and odd $\lambda$, if $\lambda \leq d + 1$, then equation (18) has no solution under (20). If $\lambda \geq d + 3$, then equation (18) has at least one solution under (20) and the complete set $S_2$ of solutions is given by*

$$k \in Z[0, (\lambda - d - 3)/4],$$
$$u_m \in O[-(\lambda - 4k - 2)/(m + 1), (\lambda - 4k - 2)/(m - 1)],$$
$$u_{m+1} \in O[(m - 1)|u_m| - u_m - (\lambda - 4k), -|u_m| - 2],$$

TABLE 1
*The number $f(n)$ of nonisomorphic $OA(n,4,2,2)$'s when $n=4\lambda$ is not a multiple of 8. Note that $f(n)=g(\lambda,2)$*

| $n$ | 12 | 20 | 28 | 36 | 44 | 52 | 60 | 68 | 76 |
|---|---|---|---|---|---|---|---|---|---|
| $f(n)$ | 1 | 3 | 7 | 15 | 28 | 48 | 79 | 123 | 184 |
| $n$ | 84 | 92 | 100 | 108 | 116 | 124 | 132 | 140 | 148 |
| $f(n)$ | 268 | 379 | 523 | 709 | 943 | 1234 | 1594 | 2032 | 2560 |
| $n$ | 156 | 164 | 172 | 180 | 188 | 196 | 204 | | |
| $f(n)$ | 3194 | 3946 | 4832 | 5872 | 7082 | 8482 | 10097 | | |

$$u_{m-1} \in O[-(\lambda - 4k + u_m + u_{m+1})/(m-1), -|u_m|],$$
$$u_j \in O[-(\lambda - 4k + u_{j+1} + \cdots + u_{m+1})/j, u_{j+1}]$$
$$\text{for } j = m-2, m-3, \ldots, 2,$$
$$u_1 = -(\lambda - 4k + u_2 + \cdots + u_{m+1}).$$

The thorough proofs of Lemmas 7 and 8 involve some tedious, yet rather straightforward, algebra and are thus omitted here. These two lemmas, in association with Theorem 1, allow the complete set of nonisomorphic $OA(\lambda 2^d, d+2, 2, d)$'s with even $d$ and odd $\lambda$ to be explicitly and systematically obtained through their $\mathbf{J}^*$-vectors given by $J_{t_j} = 2^d u_j$ for $j = 1, \ldots, m+1$, where $(u_1, \ldots, u_{m+1}, k) \in S = S_1 \cup S_2$, with $S_1$ and $S_2$ given in Lemmas 7 and 8, respectively.

Before moving on, we present in Table 1 an application of Theorem 1 and Lemmas 7 and 8 to orthogonal arrays of strength 2. Due to spatial limitations, we only present the number of nonisomorphic orthogonal arrays for $12 \leq n = 4\lambda \leq 204$ with $\lambda = 3, 5, 7, \ldots, 51$.

Some results on orthogonal arrays of strength 4 are presented in Table 2. The two cases involving $OA(80, 6, 2, 4)$'s and $OA(112, 6, 2, 4)$'s have also been solved in a recent paper by Bulutoglu and Margot [2]. Their results on $OA(144, m, 2, 4)$'s are that there are precisely 20 nonisomorphic $OA(144, 8, 2, 4)$'s and there exists no $OA(144, 9, 2, 4)$.

TABLE 2
*The number $f(n)$ of nonisomorphic $OA(n,6,2,4)$'s when $n=16\lambda$ is not a multiple of 32. Note that $f(n)=g(\lambda,4)$*

| $n$ | 80 | 112 | 144 | 176 | 208 | 240 | 272 | 304 | 336 |
|---|---|---|---|---|---|---|---|---|---|
| $f(n)$ | 1 | 3 | 7 | 14 | 26 | 46 | 77 | 123 | 190 |
| $n$ | 368 | 400 | 432 | 464 | 496 | 528 | 560 | 592 | 624 |
| $f(n)$ | 285 | 418 | 599 | 842 | 1163 | 1582 | 2123 | 2813 | 3684 |

3.2. *Enumerating $OA(\lambda 2^d, d+2, 2, d)$'s for even $d$ and even $\lambda$.* Let $\lambda = 2\lambda^*$. Since $\lambda$ is even, $\lambda^*$ is an integer. Now, consider the equation

$$\lambda^* + u_1 + \cdots + u_{m+1} = 2k, \tag{22}$$

where $k \geq 0$ is an integer and the $u_j$'s are integers with $|u_j| \leq \lambda^*$ for all $j = 1, \ldots, m+1$. Further, the $u_j$'s must satisfy either (23) or (24):

$$u_1 \leq \cdots \leq u_m \leq -|u_{m+1}|, \tag{23}$$

$$u_1 \leq \cdots \leq u_{m-1} \leq -|u_m|, \qquad u_{m+1} \leq -|u_m| - 1. \tag{24}$$

THEOREM 2. *Suppose that $d$ is even and $\lambda = 2\lambda^*$ is also even. We then have that:*

(i) *each solution $(u_1, \ldots, u_{m+1}, k)$ to equation (22) under either (23) or (24) determines an $OA(\lambda 2^d, d+2, 2, d)$ with $\mathbf{J}^*$-vector given by $J_{t_j} = 2^{d+1} u_j$ for $j = 1, \ldots, m+1$;*

(ii) *the complete set of nonisomorphic $OA(\lambda 2^d, d+2, 2, d)$'s for even $d$ and even $\lambda$ is given by collecting the arrays obtained in (i) over all the solutions to equation (22).*

We note that, unlike Theorem 1, $J_{t_j}$ and $u_j$ are now related through $J_{t_j} = 2^{d+1} u_j$. The proof of Theorem 2 is quite similar to that of Theorem 1. Here we only give an outline of the proof. To prove part (i) of Theorem 2, for $J_{t_j} = 2^{d+1} u_j$ we need to show that

$$N_s = 2^{-m}\left(n + \sum_{j=1}^{m+1} h_{st_j} J_{t_j}\right) = \left(\lambda^* + \sum_{j=1}^{m+1} h_{st_j} u_j\right)\bigg/2 = k - \sum_{j: h_{st_j} = -1} u_j$$

is a nonnegative integer for all $s \subseteq Z_m$, which must be true because of Lemma 6 and (23) or (24). Part (ii) of Theorem 2 follows from Lemma 3(ii) and Proposition 1.

EXAMPLE 2. Consider $OA(8\lambda^*, 4, 2, 2)$'s. For $\lambda^* = 1$, equation (22) has one solution given by $(u_1, \ldots, u_5, k) = (-1, 0, 0, 0, 0, 0)$ under (23) and another solution $(0, 0, 0, 0, -1, 0)$ under (24). For $\lambda^* = 2$, equation (22) has three solutions $(0, 0, 0, 0, 0, 1)$, $(-1, -1, 0, 0, 0, 0)$, $(-2, 0, 0, 0, 0, 0)$ under (23) and two solutions $(-1, 0, 0, 0, -1, 0)$, $(0, 0, 0, 0, -2, 0)$ under (24). When $\lambda^* = 3$, equation (22) has the following ten solutions under either (23) or (24):

$$(-1, 0, 0, 0, 0, 1), (0, 0, 0, 0, -1, 1), (-1, -1, -1, 0, 0, 0), (-1, -1, 0, 0, -1, 0),$$

$$(-2, -1, 0, 0, 0, 0), (-2, 0, 0, 0, -1, 0), (-1, 0, 0, 0, -2, 0), (-3, 0, 0, 0, 0, 0),$$

$$(0, 0, 0, 0, -3, 0), (-1, -1, -1, -1, 1, 0).$$

Similarly, the number of nonisomorphic $OA(32\lambda^*, 6, 2, 4)$'s is 2 and 5 for $\lambda^* = 1$ and 2, respectively. For $\lambda^* = 3$, the number of nonisomorphic $OA(32\lambda^*, 6, 2, 4)$'s is nine, with the corresponding solutions to (22) given by

$$(-1, 0, 0, 0, 0, 0, 0, 1), (0, 0, 0, 0, 0, 0, -1, 1), (-1, -1, -1, 0, 0, 0, 0, 0),$$

$$(-1, -1, 0, 0, 0, 0, -1, 0), (-2, -1, 0, 0, 0, 0, 0, 0), (-2, 0, 0, 0, 0, 0, -1, 0),$$

$$(-1, 0, 0, 0, 0, 0, -2, 0), (-3, 0, 0, 0, 0, 0, 0, 0), (0, 0, 0, 0, 0, 0, -3, 0).$$

Except for the two extra zeros in the fifth and sixth positions, these nine solutions match the first nine solutions for $OA(8\lambda^*, 4, 2, 2)$'s for $\lambda^* = 3$. These results generalize to the situation of any $d \geq 4$.

COROLLARY 3. *Let $g(\lambda, d)$ be the number of nonisomorphic $OA(\lambda 2^d, d+2, 2, d)$'s with even $d$ and even $\lambda$. Then $g(\lambda, d) \geq 1$ for all even $\lambda$. When $d = 2$, we have $g(\lambda, d) = 2, 5, 10$ for $\lambda = 2, 4, 6$, respectively. When $d \geq 4$, we have $g(\lambda, d) = 2, 5, 9$ for $\lambda = 2, 4, 6$, respectively.*

One can verify Corollary 3 directly from Theorem 2, or derive it using the following Lemmas 9 and 10, which provide all solutions to equation (22) and are the respective counterparts of Lemmas 7 and 8 in Section 3.1. Let $S_3$ denote the set of solutions to (22) under (23), and $S_4$ the set of solutions to (22) under (24).

LEMMA 9. *For even $d$ and even $\lambda$, the complete set $S_3$ of solutions to (22) under (23) is given by*

$$k \in Z[0, \lambda^*/2],$$
$$u_{m+1} \in Z[-(\lambda^* - 2k)/(m+1), (\lambda^* - 2k)/(m-1)],$$
$$u_m \in Z[-(\lambda^* - 2k + u_{m+1})/m, -|u_{m+1}|],$$
$$u_j \in Z[-(\lambda^* - 2k + u_{j+1} + \cdots + u_{m+1})/j, u_{j+1}]$$
$$\text{for } j = m-1, m-2, \ldots, 2,$$
$$u_1 = -(\lambda^* - 2k + u_2 + \cdots + u_{m+1}).$$

LEMMA 10. *For even $d$ and even $\lambda$, the complete set $S_4$ of solutions to (22) under (24) is given by*

$$k \in Z[0, (\lambda^* - 1)/2],$$
$$u_m \in Z[-(\lambda^* - 2k - 1)/(m+1), (\lambda^* - 2k - 1)/(m-1)],$$
$$u_{m+1} \in Z[(m-1)|u_m| - u_m - (\lambda^* - 2k), -|u_m| - 1],$$

$$u_{m-1} \in Z[-(\lambda^* - 2k + u_m + u_{m+1})/(m-1), -|u_m|],$$
$$u_j \in Z[-(\lambda^* - 2k + u_{j+1} + \cdots + u_{m+1})/j, u_{j+1}]$$
$$\text{for } j = m-2, m-3, \ldots, 2,$$
$$u_1 = -(\lambda^* - 2k + u_2 + \cdots + u_{m+1}).$$

Lemmas 9 and 10, together with Theorem 2, allow all nonisomorphic $OA(\lambda 2^d, d+2, 2, d)$'s with even $d$ and even $\lambda$ to be explicitly and systematically constructed through their $\mathbf{J}^*$-vectors given by $J_{t_j} = 2^{d+1} u_j$ for $j = 1, \ldots, m+1$, where $(u_1, \ldots, u_{m+1}, k) \in S_3 \cup S_4$ with $S_3$ and $S_4$ given in Lemmas 9 and 10, respectively.

Table 3 presents the number of nonisomorphic $OA(n, 4, 2, 2)$'s for $8 \leq n = 8\lambda^* \leq 200$. For each $n$, some of the arrays in Table 3 have strength 3. According to a result of Seiden and Zemach ([17]; see also [13], page 35), the number of nonisomorphic $OA(n, 4, 2, 3)$'s is precisely $[\lambda^*/2] + 1$, where $n = 8\lambda^*$. If $n$ is a multiple of 16, then exactly one of these strength 3 arrays also has strength 4. For $8 \leq n \leq 100$, Li, Deng and Tang [15] identified the nonisomorphic $OA(n, 4, 2, 2)$'s with distinct confounding frequency vectors (CFVs). Comparing their Table 1 with our Tables 1 and 3, we see that the CFV does a superb job for the case when $n$ is not a multiple of 8, as it is capable of discriminating all the $OA(n, 4, 2, 2)$'s. However, the CFV is not so discriminating when $n$ is a multiple of 8. For example, there are two nonisomorphic $OA(40, 4, 2, 2)$'s that have the same CFV, these two arrays being given by the solutions to (22),

$$(u_1, \ldots, u_5, k) = (-1, -1, -1, -1, -1, 0), (-1, -1, -1, -1, 1, 1).$$

Some results on orthogonal arrays of strength 4 are given in Table 4. Bulutoglu and Margot [2] found that the number of nonisomorphic $OA(96, m, 2, 4)$'s is 4 and 0 for $m = 7$ and 8, respectively.

TABLE 3
*The number $f(n)$ of nonisomorphic $OA(n, 4, 2, 2)$'s when $n = 4\lambda$ is a multiple of 8. Note that $f(n) = g(\lambda, 2)$*

| $n$    | 8    | 16   | 24   | 32   | 40   | 48   | 56   | 64   | 72   |
|--------|------|------|------|------|------|------|------|------|------|
| $f(n)$ | 2    | 5    | 10   | 19   | 32   | 54   | 84   | 128  | 188  |
| $n$    | 80   | 88   | 96   | 104  | 112  | 120  | 128  | 136  | 144  |
| $f(n)$ | 270  | 376  | 517  | 694  | 919  | 1198 | 1543 | 1960 | 2468 |
| $n$    | 152  | 160  | 168  | 176  | 184  | 192  | 200  |      |      |
| $f(n)$ | 3072 | 3792 | 4640 | 5636 | 6792 | 8137 | 9682 |      |      |

3.3. *Enumerating $OA(\lambda 2^d, d+2, 2, d)$'s for odd $d$ and odd $\lambda$.* Consider the equation

(25) $$\lambda + u_1 + \cdots + u_{m+1} = 4k,$$

where $k$ and the $u_j$'s are unknown integers, $k \geq 0$, $u_1, \ldots, u_m$ are odd with $|u_j| \leq \lambda - 2$ and $u_{m+1}$ is even with $|u_{m+1}| \leq \lambda - 1$. Further, the $u_j$'s must satisfy

(26) $$u_1 \leq \cdots \leq u_{m-1} \leq -|u_m|, \qquad u_{m+1} \leq 0.$$

THEOREM 3. *Consider $OA(\lambda 2^d, d+2, 2, d)$'s for odd $d$ and odd $\lambda$. We then have that:*

(i) *each solution $(u_1, \ldots, u_{m+1}, k)$ to equation (25) under (26) determines an $OA(\lambda 2^d, d+2, 2, d)$ with $\mathbf{J}^*$-vector given by $J_{t_j} = 2^d u_j$ for $j = 1, \ldots, m+1$;*

(ii) *the complete set of nonisomorphic $OA(\lambda 2^d, d+2, 2, d)$'s for odd $d$ and odd $\lambda$ is given by collecting the arrays obtained in (i) over all the solutions to equation (25) under (26).*

We omit the proof of Theorem 3, which is similar to that of Theorem 1 and involves the use of Lemma 1, Lemma 3(iv), Lemma 4, Lemma 6 and Proposition 2.

EXAMPLE 3. Consider $OA(8\lambda, 5, 2, 3)$'s. For $\lambda = 1$, equation (25) has no solution under (26). For $\lambda = 3$, equation (25) has a unique solution under (26) given by $(u_1, \ldots, u_6, k) = (-1, -1, -1, -1, 1, 0, 0)$. When $\lambda = 5$, there are three solutions $(-3, -1, -1, -1, 1, 0, 0)$, $(-1, -1, -1, -1, -1, 0, 0)$ and $(-1, -1, -1, -1, 1, -2, 0)$. When $\lambda = 7$, there are seven solutions:

$$(-1, -1, -1, -1, 1, 0, 1), (-5, -1, -1, -1, 1, 0, 0), (-1, -1, -1, -1, 1, -4, 0),$$

$$(-3, -3, -1, -1, 1, 0, 0), (-3, -1, -1, -1, -1, 0, 0),$$

$$(-3, -1, -1, -1, 1, -2, 0), (-1, -1, -1, -1, -1, -2, 0).$$

These generalize to $OA(\lambda 2^d, d+2, 2, d)$'s with $\lambda = d-2, d, d+2, d+4$ for any $d \geq 3$.

TABLE 4
*The number $f(n)$ of nonisomorphic $OA(n, 6, 2, 4)$'s when $n = 16\lambda$ is a multiple of 32. Note that $f(n) = g(\lambda, 4)$*

| $n$ | 32 | 64 | 96 | 128 | 160 | 192 | 224 | 256 | 288 |
|---|---|---|---|---|---|---|---|---|---|
| $f(n)$ | 2 | 5 | 9 | 17 | 29 | 49 | 77 | 120 | 179 |
| $n$ | 320 | 352 | 384 | 416 | 448 | 480 | 512 | 544 | 576 |
| $f(n)$ | 265 | 380 | 539 | 747 | 1025 | 1383 | 1848 | 2435 | 3181 |

TABLE 5
*The number $f(n)$ of nonisomorphic $OA(n, 5, 2, 3)$'s when $n = 8\lambda$ is not a multiple of 16. Note that $f(n) = g(\lambda, 3)$*

| $n$ | 24 | 40 | 56 | 72 | 88 | 104 | 120 | 136 | 152 |
|---|---|---|---|---|---|---|---|---|---|
| $f(n)$ | 1 | 3 | 7 | 15 | 28 | 49 | 82 | 130 | 199 |
| $n$ | 168 | 184 | 200 | 216 | 232 | 248 | 264 | 280 | 296 |
| $f(n)$ | 296 | 428 | 605 | 839 | 1142 | 1530 | 2022 | 2637 | 3399 |
| $n$ | 312 | 328 | 344 | 360 | 376 | 392 | 408 | | |
| $f(n)$ | 4336 | 5476 | 6854 | 8509 | 10481 | 12818 | 15573 | | |

COROLLARY 4. *Let $g(\lambda, d)$ be the number of nonisomorphic $OA(\lambda 2^d, d+2, 2, d)$'s with odd $d$ and odd $\lambda$. Then $g(\lambda, d) = 0$ for $\lambda \leq d - 2$ and $g(\lambda, d) = 1, 3, 7$ for $\lambda = d, d+2, d+4$, respectively.*

Corollary 4 can be verified directly using Theorem 3, or derived from the following Lemma 11, which gives the complete set of solutions to (25) under (26). To present Lemma 11, we need to introduce the notation $E[a, b]$, which represents the set of even integers $x$ satisfying $a \leq x \leq b$.

LEMMA 11. *For odd $d$ and odd $\lambda$, if $\lambda \leq d - 2$, then equation (25) has no solution under (26). If $\lambda \geq d$, then equation (25) has at least one solution under (26) and the complete set of solutions is given by*

$$k \in Z[0, (\lambda - d)/4],$$

$$u_{m+1} \in E[-(\lambda - d - 4k), 0],$$

$$u_m \in O[-(\lambda - 4k + u_{m+1})/m, (\lambda - 4k + u_{m+1})/(m-2)],$$

$$u_{m-1} \in O[-(\lambda - 4k + u_m + u_{m+1})/(m-1), -|u_m|],$$

$$u_j \in O[-(\lambda - 4k + u_{j+1} + \cdots + u_{m+1})/j, u_{j+1}]$$

$$\text{for } j = m - 2, m - 2, \ldots, 2,$$

$$u_1 = -(\lambda - 4k + u_2 + \cdots + u_{m+1}).$$

Based on Theorem 3 and Lemma 11, all nonisomorphic $OA(\lambda 2^d, d+2, 2, d)$'s for given odd $d$ and odd $\lambda$ can be obtained. In Table 5, we present the number of nonisomorphic $OA(8\lambda, 5, 2, 3)$'s for $\lambda = 3, 5, \ldots, 51$.

3.4. *Enumerating $OA(\lambda 2^d, d+2, 2, d)$'s for odd $d$ and even $\lambda$.* Let $\lambda = 2\lambda^*$. Since $\lambda$ is even, $\lambda^*$ is an integer. Now consider the equation

(27) $$\lambda^* + u_1 + \cdots + u_{m+1} = 2k,$$

where $k$ and the $u_j$'s are unknown integers, $k \geq 0$ and $|u_j| \leq \lambda^*$ for $j = 1, \ldots, m+1$. Further, the $u_j$'s must satisfy

(28) $$u_1 \leq \cdots \leq u_{m-1} \leq -|u_m|, \qquad u_{m+1} \leq 0.$$

THEOREM 4. *Consider $OA(\lambda 2^d, d+2, 2, d)$'s for odd $d$ and even $\lambda = 2\lambda^*$. We then have that:*

(i) *each solution $(u_1, \ldots, u_{m+1}, k)$ to equation (27) under (28) determines an $OA(\lambda 2^d, d+2, 2, d)$ with $\mathbf{J}^*$-vector given by $J_{t_j} = 2^{d+1} u_j$ for $j = 1, \ldots, m+1$;*

(ii) *the complete set of nonisomorphic $OA(\lambda 2^d, d+2, 2, d)$'s for odd $d$ and even $\lambda$ is given by collecting the arrays obtained in (i) over all the solutions to equation (27) under (28).*

Theorem 4 is similar to Theorem 2 in that $J_{t_j}$ and $u_j$ are related through $J_{t_j} = 2^{d+1} u_j$. Its proof is also straightforward. Similar to Corollary 3, the following result holds.

COROLLARY 5. *Let $g(\lambda, d)$ be the number of nonisomorphic $OA(\lambda 2^d, d+2, 2, d)$'s with odd $d$ and even $\lambda$. Then $g(\lambda, d) \geq 1$ for all even $\lambda$. When $d = 3$, we have $g(\lambda, d) = 2, 5, 10$ for $\lambda = 2, 4, 6$, respectively. When $d \geq 5$ we have $g(\lambda, d) = 2, 5, 9$ for $\lambda = 2, 4, 6$, respectively.*

All nonisomorphic $OA(\lambda 2^d, d+2, 2, d)$'s for odd $d$ and even $\lambda = 2\lambda^*$ can be obtained from Theorem 4 and the following Lemma 12, which gives all solutions to (27) under (28).

LEMMA 12. *For odd $d$ and even $\lambda$, the complete set of solutions to (27) under (28) is given by*

$$k \in Z[0, \lambda^*/2],$$
$$u_{m+1} \in Z[-(\lambda^* - 2k), 0],$$
$$u_m \in Z[-(\lambda^* - 2k + u_{m+1})/m, (\lambda^* - 2k + u_{m+1})/(m-2)],$$
$$u_{m-1} \in Z[-(\lambda^* - 2k + u_m + u_{m+1})/(m-1), -|u_m|],$$
$$u_j \in Z[-(\lambda^* - 2k + u_{j+1} + \cdots + u_{m+1})/j, u_{j+1}]$$
$$\text{for } j = m-2, m-2, \ldots, 2,$$
$$u_1 = -(\lambda^* - 2k + u_2 + \cdots + u_{m+1}).$$

Applying Theorem 4 and Lemma 12 to the case $d = 3$, we have obtained all $OA(16\lambda^*, 5, 2, 3)$'s with $\lambda^* = 1, \ldots, 25$. In Table 6, we present the number of nonisomorphic $OA(16\lambda^*, 5, 2, 3)$'s for these values of $\lambda^*$. Some of the arrays in Table 6 have strength 4 and the number of such nonisomorphic arrays is $[\lambda^*/2] + 1$ ([13], page 35). When $n$ is a multiple of 32, exactly one of these strength 4 arrays has strength 5.

**4. Discussion.** Using the theory of $J$-characteristics, we have completely solved the problem of enumerating two-level orthogonal arrays of strength $d$ with $d+2$ constraints. This further demonstrates that $J$-characteristics provide a powerful tool for studying two-level orthogonal arrays. Previously, $J$-characteristics have been successfully used in defining generalizations of minimum aberration and in providing statistical justifications for the criteria of generalized aberration.

It is tempting to consider the problem of enumerating two-level orthogonal arrays of strength $d$ with $m \geq d+3$ constraints using the approach in this paper. This seems quite nontrivial. Although the general ideas continue to apply, useful results like Theorems 1–4 may not exist after all. The following discussion should shed some light on the issue. First, let us consider enumerating orthogonal arrays of strength $d$ with $m = d+1$ constraints. Based on Lemma 2, we only need to consider the single $J$-characteristic for all $m = d+1$ columns. It is also easy to see that we can assume that this $J$ value is nonpositive. Using Lemmas 1 and 4, we see that the complete set of nonisomorphic $OA(\lambda 2^d, d+1, 2, d)$'s is given by the set of solutions to the equation $\lambda + u = 2k$, where $k \geq 0$ is an integer and $-\lambda \leq u \leq 0$. This provides a very simple derivation for the result of [17] on nonisomorphic $OA(\lambda 2^d, d+1, 2, d)$'s. As much as we like this simple derivation, it also points to the difficulty arising from trying to solve the problem with $d+3$ constraints. Judging from the huge jump, in terms of difficulty, from the case of $d+1$ to the case of $d+2$, one can only imagine what it will take to solve the case of $d+3$ completely. In our enumeration of $OA(\lambda 2^d, d+2, 2, d)$'s, Lemma 6 plays an instrumental role. It achieves two goals, one being that it substantially narrows down the possible configurations of $J$-characteristics and the other being that we only need to check that $N_s$ is a nonnegative integer for $s = \phi$. Although one can derive a result like Lemma 6 for the case of $m = d + 3$, it will not be as powerful as Lemma 6 in achieving the two aforementioned goals. Having stated all of the negatives, we recall that everything has its positive side. Completely solving the enumeration problem for the case of $m = d + 3$ may just be too ambitious after all. It is, however, possible to find a partial solution and this is what we will focus on in our future research on this problem.

TABLE 6
*The number $f(n)$ of nonisomorphic $OA(n, 5, 2, 3)$'s when $n = 8\lambda$ is a multiple of* 16.
*Note that $f(n) = g(\lambda, 3)$*

| $n$    | 16   | 32   | 48   | 64   | 80    | 96    | 112   | 128  | 144  |
|--------|------|------|------|------|-------|-------|-------|------|------|
| $f(n)$ | 2    | 5    | 10   | 19   | 33    | 56    | 89    | 138  | 207  |
| $n$    | 160  | 176  | 192  | 208  | 224   | 240   | 256   | 272  | 288  |
| $f(n)$ | 303  | 432  | 606  | 832  | 1126  | 1501  | 1975  | 2566 | 3300 |
| $n$    | 304  | 320  | 336  | 352  | 368   | 384   | 400   |      |      |
| $f(n)$ | 4198 | 5293 | 6615 | 8202 | 10092 | 12335 | 14975 |      |      |

**Acknowledgments.** The authors wish to thank an Associate Editor and two referees for their helpful comments.

DEPARTMENT OF STATISTICS  
UNIVERSITY OF GEORGIA  
ATHENS, GEORGIA 30602-1952  
USA  
E-MAIL: jstufken@uga.edu

DEPARTMENT OF STATISTICS AND ACTUARIAL SCIENCE  
SIMON FRASER UNIVERSITY  
BURNABY, BRITISH COLUMBIA  
CANADA V5A 1S6  
E-MAIL: boxint@cs.sfu.ca